\documentclass[11pt]{article}
\usepackage{amssymb, amsthm, amsmath}

\newtheorem{thm}{Theorem}[section]
\newtheorem{lem}[thm]{Lemma}
\newtheorem{prop}[thm]{Proposition}
\newtheorem{cor}[thm]{Corollary}
\theoremstyle{definition}\newtheorem{df}[thm]{Definition}
\theoremstyle{definition}\newtheorem{exm}[thm]{Example}

\renewcommand{\phi}{\varphi}

\newcommand{\N}{\mathbb{N}}
\newcommand{\Z}{\mathbb{Z}}
\newcommand{\R}{\mathbb{R}}

\newcommand{\Homeo}{\operatorname{Homeo}}
\newcommand{\md}{\operatorname{mod}}
\newcommand{\sgn}{\operatorname{sgn}}
\newcommand{\id}{\operatorname{id}}

\newcommand{\xp}{(X,\phi )}
\newcommand{\yp}{(Y,\psi )}

\title{Some remarks on topological full groups \\ 
of Cantor minimal systems}
\author{MATUI Hiroki
\thanks{The author was supported by 
Grant-in-Aid for Young Scientists (B) of
Japan Society for the Promotion of Science.}}
\date{}

\begin{document}
\maketitle
\begin{abstract}
Giordano, Putnam and Skau showed that topological full groups 
of Cantor minimal systems are complete invariants for flip conjugacy. 
We will completely determine the structure of normal subgroups 
of the topological full group. 
Moreover, a necessary and sufficient condition for the topological 
full group to be finitely generated will be given. 
\end{abstract}

\section{Introduction}
The notion of full group was introduced by Henry Dye 
in the measurable dynamical setting. 
He showed that if $G_1$ and $G_2$ are two countable groups of 
measure preserving transformations acting ergodically 
on a Lebesgue space, any group isomorphism between the full groups 
$[G_1]$ and $[G_2]$ is implemented by an orbit equivalence of 
$G_1$ and $G_2$. 
Giordano, Putnam and Skau studied its topological analogue in 
\cite{GPS2}. 
They introduced several kinds of full groups 
for Cantor minimal systems and showed that they are complete 
invariants for orbit equivalence, strong orbit equivalence and 
flip conjugacy. 
Here, by a Cantor minimal system, we mean a pair of a Cantor set 
$X$ and a minimal homeomorphism $\phi\in\Homeo(X)$. 
The Cantor set $X$ is characterized as a topological space 
which is compact, metrizable, perfect and totally disconnected. 
A homeomorphism $\phi$ is said to be minimal, 
if there is no non-trivial closed $\phi$-invariant closed subset. 

The purpose of this paper is to clarify the structure of 
topological full groups associated with Cantor minimal systems. 
The topological full group $[[\phi]]$ consists of all homeomorphisms 
$\gamma$ that are locally equal to some power of 
the minimal homeomorphism $\phi$ (see Definition \ref{fg} below). 
As the number of clopen sets in $X$ is only countable, 
the set of $\Z$-valued continuous functions on $X$ 
denoted by $C(X,\Z)$ is also countable. 
Thus $[[\phi]]$ is a countable group. 
This is a great difference from $[\phi]$: 
the full group $[\phi]$ is always uncountable, 
and a complete invariant for orbit equivalence, 
although $[[\phi]]$ is a complete invariant for flip conjugacy. 
Since the topological full group $[[\phi]]$ has full information of 
$\xp$ and is also countable, 
it is natural and important to investigate its various properties 
as a group. 
The present paper is the first attempt at such a study. 

Now we give an overview of each section below. 
In Section 2, we will recall several results from $\cite{GPS1}$, 
$\cite{GPS2}$ and $\cite{GPS3}$. 
We write the kernel of the index map $I:[[\phi]]\to\Z$ by $[[\phi]]_0$ 
(see Definition \ref{index} below). 
In Section 3, the notion of AF full groups is introduced and 
the structure of its normal subgroups is discussed. 
Based on the results obtained in Section 3, 
we will define the signature on $[[\phi]]_0$ in Section 4. 
The signature takes its value in $K^0\xp/2K^0\xp$, 
where $K^0\xp$ is defined as 
\[ K^0\xp \ = \ C(X,\Z)/\{f-f\phi^{-1}:f\in C(X,\Z)\}. \]
Let $[f]$ denote the equivalence class of a function $f$ 
in $K^0\xp$. 
This $K^0$-group played a central role in \cite{GPS1}. 
They showed that $K^0$-groups of two Cantor minimal systems are 
unital order isomorphic if and only if the two systems are 
strong orbit equivalent. 
In the present paper we show that the quotient of $[[\phi]]_0$ 
by its commutator subgroup $D([[\phi]]_0)$ is isomorphic to 
the $\Z_2$-valued dynamical cohomology $K^0\xp\otimes\Z_2
\cong K^0\xp/2K^0\xp$. 
In addition, the commutator subgroup $D([[\phi]]_0)$ is shown 
to be simple, which was mentioned in \cite[Theorem 4.11]{S2} 
without proof. 
In Section 5, it is proved that $D([[\phi]]_0)$ is finitely 
generated if and only if $\xp$ is a minimal subshift. 
We need the result of Section 4 to determine when $[[\phi]]$ is 
finitely generated. 
It is also shown that $D([[\phi]]_0)$ never be finitely presented. 
In Section 6, two examples of minimal subshifts are given. 

\section{Preliminaries}
We would like to collect basic facts which are needed 
in later sections. 
We will consider a countable equivalence relation ${\cal R}$ 
on a compact metrizable zero-dimensional space $X$, that is, 
$X$ is homeomorphic to a compact subset of $\{0,1\}^\N$ with 
the product topology and ${\cal R}\subset X\times X$ is 
an equivalence relation so that each equivalence class 
$[x]_{\cal R}=\{y\in X:(x,y)\in {\cal R}\}$ is countable 
for $x\in X$. 
If $(x,y),(y,z)\in {\cal R}$, 
then these pairs are said to be composable and their product is 
defined as $(x,y)(y,z)=(x,z)$. For $(x,y)\in {\cal R}$ its inverse 
is defined as $(x,y)^{-1}=(y,x)$. With this structure ${\cal R}$ is 
a (principal) groupoid. We say that ${\cal R}$ is minimal, if 
$[x]_{\cal R}$ is dense in $X$ for all $x\in X$. 
\begin{df}[{\cite[Definition 2.1, 3.1, 3.7]{GPS3}}]
Let $X$ be a compact metrizable zero-dimensional space and 
let ${\cal R}\subset X\times X$ be a countable equivalence relation. 
\begin{enumerate}
\item Suppose that ${\cal R}$ is equipped 
with a Hausdorff locally compact, 
second countable topology ${\cal T}$ so that the product of composable 
pairs is continuous and and the inverse map is a homeomorphism. 
We say that $({\cal R},{\cal T})$ is an \'{e}tale equivalence 
relation on $X$, if the range map $r$ defined by 
${\cal R}\ni (x,y)\mapsto x\in X$ is a local homeomorphism. 
\item An \'{e}tale equivalence relation $({\cal R},{\cal T})$ on $X$ 
is called a compact \'{e}tale equivalence relation (CEER for short), 
if ${\cal R}$ is compact. 
\item Let $\{({\cal R}_n,{\cal T}_n)\}_{n=1}^\infty$ be a sequence 
of CEER's on $X$, so that ${\cal R}_n\subset {\cal R}_{n+1}$ and 
${\cal R}_n\in {\cal T}_{n+1}$. Let $({\cal R}, {\cal T})$ be 
the inductive limit of $\{({\cal R}_n,{\cal T}_n)\}_{n=1}^\infty$ 
with the inductive limit topology ${\cal T}$. 
We say that $({\cal R}, {\cal T})$ is an AF equivalence relation and 
use the notation $({\cal R}, {\cal T})=\lim ({\cal R}_n,{\cal T}_n)$. 
\end{enumerate}
We remark that if $({\cal R}, {\cal T})$ is a CEER, ${\cal T}$ is 
the relative topology from $X\times X$ and 
the cardinality of each equivalence class is uniformly finite. 
\end{df}
\begin{df}[{\cite[Example 2.7(ii)]{GPS3}}]
We say that $B=(V,E)$ is a Bratteli diagram when $V=\bigcup _{n=0}
^\infty V_n$ and $E=\bigcup _{n=1}^\infty E_n$ are disjoint unions 
of finite sets of vertices and edges with source maps 
$s:E_n\rightarrow V_{n-1}$ and range maps $r:E_n\rightarrow V_n$ 
both of which are surjective. 
We always assume that $V_0$ consists of one 
point $v_0$. For a Bratteli diagram $B=(V,E)$, 
\[ X_B=\{(e_n)_n:e_n\in E_n, r(e_n)=s(e_{n+1}) \ 
\mbox{for all} \ n\in \N \}\]
is called the infinite path space of $B$. The space $X_B$ is 
endowed with the natural product topology. 

We say that $B=(V,E)$ is simple when 
for every $v\in V_n$ there exists $m>n$ such that all vertices 
of $V_m$ are connected to $v$. 
\end{df}
When $B=(V,E)$ is a Bratteli diagram, the infinite path space $X_B$ 
is a compact metrizable zero-dimensional space. 
Observe that if $B$ is simple, $X_B$ is a finite set or a Cantor set. 
\begin{df}[{\cite[Example 2.7(ii)]{GPS3}}]
Let $B=(V,E)$ be a Bratteli diagram. 
Two infinite paths are said to be tail equivalent, 
if they agree from some level on. 
We denote the tail equivalence relation by ${\cal R}_B$. 
Clearly ${\cal R}_B$ is an increasing union of 
\[ {\cal R}_N \ = \ \{((e_n)_n,(f_n)_n)\in X_B\times X_B:
e_m=f_m \mbox{ for all } m\geq N\}\]
for $N\in \N$. Given the relative topology of $X\times X$, 
${\cal R}_N$ is a CEER. 
Give ${\cal R}_B$ the inductive limit topology ${\cal T}_B$. 
We call $({\cal R}_B, {\cal T}_B)$ the AF equivalence relation 
associated with $B$. 
\end{df}
The following theorem asserts that every AF equivalence relation 
can be expressed by a Bratteli diagram. 
\begin{thm}[{\cite[Theorem 3.9]{GPS3}}]\label{AFmodel}
Let $({\cal R},{\cal T})$ be an AF equivalence relation on $X$. 
There exists a Bratteli diagram $B=(V,E)$ such that 
$({\cal R},{\cal T})$ is isomorphic to the AF equivalence relation 
$({\cal R}_B,{\cal T}_B)$ associated with $B$. 
Thus, there exists a homeomorphism $\pi :X\rightarrow X_B$ 
such that $\pi\times \pi$ gives a homeomorphism from 
${\cal R}$ to ${\cal R}_B$. Furthermore, ${\cal R}$ is minimal 
if and only if $B$ is simple. 
\end{thm}
Let us recall that a dimension group $K_0(B)$ is associated to 
a Bratteli diagram $B=(V,E)$. 
The notation is derived from the $K$-theory of $C^*$-algebras: 
the Bratteli diagram $B=(V,E)$ gives rise to an approximately 
finite-dimensional (AF) $C^*$-algebra and 
its $K_0$-group is isomorphic to $K_0(B)$. 
The reader may refer to Section 3 of \cite{GPS1} 
for the connection to $K$-theory of $C^*$-algebras. 
The dimension group $K_0(B)$ is defined as the inductive limit 
of a free abelian group $\Z^{V_n}$ of finite rank. 
The connecting map from $\Z^{V_n}$ to $\Z^{V_{n+1}}$ is given by 
\[ \Z^{V_n}\ni v \ \mapsto \ \sum_{w\in V_{n+1}}
m_{v,w}w\in \Z^{V_{n+1}}, \]
where the canonical basis is denoted by the vertices and 
$m_{v,w}$ is the number of edges which go from $v\in V_n$ to 
$w\in V_{n+1}$. Endowed with a positive cone $K_0(B)^+$ and 
an order unit, $K_0(B)$ is a unital ordered group. 
We do not go to the detail, because it is not needed 
in this paper. 
\bigskip

We will review relation between Cantor minimal systems and 
AF equivalence relations. 
Let $\phi$ be a minimal homeomorphism on the Cantor set $X$. 
As explained in \cite[Example 2.7(i)]{GPS3}, 
an \'{e}tale equivalence relation $({\cal R},{\cal T})$ is 
associated with $\xp$. 
More precisely, 
\[ {\cal R} \ = \ \{(x,\phi^n(x)):x\in X, \ n\in \Z\} \]
and ${\cal T}$ is given by transferring the product topology 
on $X\times \Z$ to ${\cal R}$ via the bijection 
$(x,n)\mapsto (x,\phi^n(x))$. 

A subset $Y\subset X$ is said to be wandering 
if $Y\cap \phi^n(Y)$ is empty for all $n\in \N$. 
We consider the subequivalence relation ${\cal R}_Y$ of ${\cal R}$ 
defined by 
\[ {\cal R}_Y \ = \ \{(x,\phi^n(x)),(\phi^n(x),x):
x\in X, \ n\geq 0, \ \phi^k(x)\notin Y \ 
\mbox{ if }0\leq k< n-1\} \]
for a closed wandering subset $Y$. 
Thus ${\cal R}_Y$ is obtained from ${\cal R}$ by splitting 
the $\phi$-orbits meeting $Y$ in the forward and backward orbits 
at $Y$. 
Let ${\cal T}_Y$ denote the relative topology from ${\cal T}$. 
Although ${\cal R}_Y$ is defined in Section 4 of \cite{GPS3} 
in the case that $Y$ is not necessarily wandering, 
we do not need that case here. 
If $Y$ is a singleton $\{x\}$, we write $({\cal R}_x, {\cal T}_x)$ 
instead of $({\cal R}_{\{x\}}, {\cal T}_{\{x\}})$. 
\begin{thm}[{\cite[Lemma 5.1]{GPS1}}{\cite[Theorem 4.6]{GPS3}}]
\label{AFinCM}
Let $\xp$ be a Cantor minimal system and let $({\cal R}, {\cal T})$ 
be the associated \'{e}tale equivalence relation. When $Y\subset X$ 
is a non-empty closed wandering set, 
the equivalence relation $({\cal R}_Y,{\cal T}_Y)$ is AF. 
\end{thm}
The key of the above theorem is that $({\cal R}_x,{\cal T}_x)$ is 
an AF equivalence relation. This fact goes back to \cite{HPS} and 
it enables us to represent a minimal homeomorphism on the Cantor 
set as a Bratteli-Vershik map on the infinite path space of 
a Bratteli diagram. 
Besides, under this representation, we can see that 
the dimension group associated with the Bratteli diagram is 
isomorphic to $K^0\xp$. We refer the reader to \cite{HPS}. 
\bigskip

In the last part of this section, we will review three kinds 
of full groups which were introduced in \cite{GPS2}. 
\begin{df}[{\cite[Definition 1.1, 2.1]{GPS2}}]\label{fg}
Let $\xp$ be a Cantor minimal system. 
\begin{enumerate}
\item The full group $[\phi]$ of $\xp$ is the subgroup of 
all homeomorphisms $\gamma\in\Homeo(X)$ such that 
$\gamma(x)\in\{\phi^n(x):n\in\Z\}$ for all $x\in X$. 
To any $\gamma\in[\phi]$ is associated a map $n:X\rightarrow \Z$, 
defined by $\gamma(x)=\phi^{n(x)}(x)$ for $x\in X$. 
\item The topological full group $[[\phi]]$ of $\xp$ is 
the subgroup of all homeomorphisms $\gamma\in\Homeo(X)$, 
whose associated map $n:X\rightarrow \Z$ is continuous. 
\item For every $x\in X$, let $[[\phi]]_x$ denote the subgroup 
of all elements $\gamma\in[[\phi]]$ that preserve the forward 
orbit $\{\phi^n(x):n\in\N\}$ and the backward orbit 
$\{\phi^{1-n}(x):n\in\N\}$ respectively. 
\end{enumerate}
\end{df}
Let $({\cal R},{\cal T})$ be the equivalence relation associated 
with a Cantor minimal system $\xp$. 
Then $[\phi]$ consists of homeomorphisms $\gamma$ which preserves 
every equivalence class $[x]_{\cal R}$, and 
$[[\phi]]$ consists of the homeomorphism $\gamma$ such that 
$X\ni x\mapsto (x,\gamma(x))\in{\cal R}$ is well-defined and 
continuous. 
For every $x\in X$, $[[\phi]]_x$ is the set of all homeomorphisms 
$\gamma\in[[\phi]]$ that preserves 
every ${\cal R}_x$-equivalence class. 
\begin{df}[{\cite[Section 5]{GPS2}}]\label{index}
Let $\xp$ be a Cantor minimal system and $\mu$ be a $\phi$-invariant 
probability measure on $X$. For $\gamma\in[[\phi]]$, we define 
\[ I(\gamma) \ = \ \int n_{\gamma} \ d\mu, \]
where $n_{\gamma}$ is the continuous map from $X$ to $\Z$ 
associated with $\gamma$. 
Then $I(\gamma)$ does not depend on the choice of $\mu$ and 
$I$ is a surjective homomorphism from $[[\phi]]$ to $\Z$. 
We call the homomorphism $I$ the index map and 
denote its kernel by $[[\phi]]_0$. 
\end{df}
The reason why the value of the index map falls into the integers 
is that the $K_1$-group of the crossed product $C^*$-algebras 
arising from $\xp$ is isomorphic to the integers. 
We refer the reader to \cite[Section 5]{GPS2} for details. 
We also remark that $[[\phi]]_x$ is contained in $[[\phi]]_0$ 
because every element of $[[\phi]]_x$ is of finite order. 

The following is the main theorem of \cite{GPS2}. 
\begin{thm}[{\cite[Corollary 4.4, 4.6, 4.11, 5.18]{GPS2}}]
\label{fulliso}
Let $\xp$ and $\yp$ be Cantor minimal systems. 
\begin{enumerate}
\item The two Cantor minimal systems are orbit equivalent 
if and only if $[\phi]$ and $[\psi]$ are isomorphic 
as abstract groups. 
\item The two Cantor minimal systems are flip conjugate 
if and only if $[[\phi]]$ and $[[\psi]]$ are isomorphic 
as abstract groups. 
\item The two Cantor minimal systems are flip conjugate 
if and only if $[[\phi]]_0$ and $[[\psi]]_0$ are isomorphic 
as abstract groups. 
\item Let $x\in X$ and $y\in Y$. 
The two Cantor minimal systems are strong orbit equivalent 
if and only if $[[\phi]]_x$ and $[[\psi]]_y$ are isomorphic 
as abstract groups. 
\end{enumerate}
\end{thm}

\section{AF full groups}
In this section we will introduce AF full groups and study their 
normal subgroups. The signature map defined on the AF full group 
will play a crucial role in the next section. 
\begin{df}
Let $X$ be a compact metrizable zero-dimensional space and let 
$({\cal R},{\cal T})$ be an \'{e}tale equivalence relation on $X$. 
We put 
\[ [{\cal R},{\cal T}] \ = \ \{ \ \gamma \in \Homeo (X):
X\ni x\mapsto (x,\gamma(x))\in {\cal R} \mbox{ is well-defined and 
continuous} \ \}\]
and call it the topological full group associated with 
$({\cal R},{\cal T})$. 
\end{df}
Since the diagonal $\{(x,x):x\in X\}$ is a clopen subset 
in ${\cal R}$, 
the fixed point set of each element $\gamma\in[{\cal R},{\cal T}]$ 
is clopen. 
When $({\cal R},{\cal T})$ is the \'{e}tale equivalence relation 
associated with a Cantor minimal system $\xp$, 
the topological full group $[{\cal R},{\cal T}]$ defined above 
equals $[[\phi]]$ defined in Section 2. 
It is also obvious that $[[\phi]]_x$ equals $[{\cal R}_x,{\cal T}_x]$ 
for all $x\in X$. 

The following proposition asserts that AF equivalence relations can 
be characterized by their associated topological full groups. 
Although the idea of the proof has already appeared 
in \cite{D} or \cite{K} 
and this fact is folklore for experts, 
we would like to give a proof without using Theorem \ref{AFmodel} 
for the reader's convenience. 
\begin{prop}
For an \'{e}tale equivalence relation $({\cal R},{\cal T})$ on the 
Cantor set $X$, the following are equivalent. 
\begin{enumerate}
\item $({\cal R},{\cal T})$ is an AF equivalence relation. 
\item $[{\cal R},{\cal T}]$ is locally finite. 
\end{enumerate}
\end{prop}
\begin{proof}
(i)$\Rightarrow$(ii) Let $({\cal R},{\cal T})$ be an AF equivalence 
relation. Suppose that $({\cal R},{\cal T})$ is the inductive limit of 
CEER's $\{({\cal R}_n,{\cal T}_n)\}_n$. 
Let $\gamma_1,\gamma_2,\dots ,\gamma_m$ be elements of 
$[{\cal R},{\cal T}]$ and 
let $H$ be the subgroup they generate. 
We must show that $H$ is finite. 
Put $F=\{\gamma_1,\gamma_1^{-1},\gamma_2,\gamma_2^{-1}, 
\dots ,\gamma_m, \gamma_m^{-1}\}$. 
Since $X$ is compact, $\{(x,\gamma_i(x)):x\in X\}$ and their union 
are also compact in ${\cal R}$. 
Hence there exists $n$ such that ${\cal R}_n$ 
contains $(x,\gamma_i(x))$ for all $x\in X$ and $i=1,2,\dots ,m$. 
Let $d$ be a metric on $X$. There is $\varepsilon>0$ such that 
$d(x,x')>2\varepsilon$ if $(x,x')\in{\cal R}_n$ and $x\neq x'$, 
because $({\cal R}_n,{\cal T}_n)$ is a CEER. 
Choose a real number $\delta >0$ smaller than $\varepsilon$ 
so that $d(x,y)<\delta$ implies $d(\gamma(x),\gamma(y))<\varepsilon$ 
for all $\gamma\in F$. 
For every $x\in X$, let $U_x$ be the open set defined by 
\[ U_x \ = \ \{y\in X:d(y',[x]_{{\cal R}_n})<\delta 
\mbox{ for all } y'\in[y]_{{\cal R}_n} \}. \]
If $y\in U_x$, there exists a unique $x'\in[x]_{{\cal R}_n}$ with 
$d(y,x')<\delta$. 
Furthermore, by the choice of $\varepsilon$, we have 
$d(\gamma(y),\gamma(x'))<\delta$ for all $\gamma\in F$, 
which implies $d(\gamma(y),\gamma(x'))<\delta$ for all $\gamma\in H$. 
Every element $\gamma$ belonging to $H$ induces a permutation 
on the finite set $[x]_{{\cal R}_n}$. 
The above argument shows that if $\gamma$ induces the identity 
on $[x]_{{\cal R}_n}$, then $\gamma$ is the identity on $U_x$. 
Since $X$ is compact, it can be covered by finitely many $U_x$'s. 
Then we obtain an injective homomorphism from $H$ 
to a finite direct sum of permutation groups, 
which completes the proof. 

(ii)$\Rightarrow$(i) Suppose that $[{\cal R},{\cal T}]$ is 
an increasing union of finite subgroups $G_1,G_2,\dots$ and 
let ${\cal R}_n=\{(x,\gamma(x)):x\in X,\gamma\in G_n\}$. 
Let ${\cal T}_n$ be the relative topology from ${\cal T}$. 
Then $({\cal R}_n,{\cal T}_n)$ is a CEER, 
because $\{(x,\gamma(x)):x\in X\}$ is compact for all 
$\gamma\in G_n$. 
To complete the proof, it suffices to show 
${\cal R}=\bigcup_n{\cal R}_n$. Let $(x,y)\in {\cal R}$. 
We may assume $x\neq y$. 
Choose a clopen neighborhood $U$ of $(x,y)$ so that the range map 
and the source map on $U$ are local homeomorphisms, and 
$r(U)$ and $s(U)$ are disjoint. Define $\gamma\in\Homeo(X)$ by 
$\gamma(z)=s_0r_0^{-1}(z)$ if $z\in r(U)$, $\gamma(z)=r_0s_0^{-1}(z)$ if 
$z\in s(U)$ and $\gamma(z)=z$ otherwise, 
where $s_0$ and $r_0$ denote the restriction to $U$ 
of $s$ and $r$, respectively. 
Then $X\ni z\mapsto (z,\gamma(z))\in {\cal R}$ is clearly 
a well-defined continuous map. 
It follows that we have $\gamma\in G_n$ for some $n$ and 
$(x,y)=(x,\gamma(x))\in {\cal R}_n$. 
\end{proof}
We call $[{\cal R},{\cal T}]$ an AF full group 
if $({\cal R},{\cal T})$ is an AF equivalence relation. 

Let $({\cal R},{\cal T})$ be an AF equivalence relation on $X$ 
and let $B=(V,E)$ be a Bratteli diagram 
such that $({\cal R},{\cal T})$ is isomorphic to the AF equivalence 
relation associated with $B$. 
We identify the infinite path space $X_B$ with $X$ and 
$({\cal R}_B,{\cal T}_B)$ with $({\cal R},{\cal T})$. 
Put $G=[{\cal R},{\cal T}]$. 
We would like to examine the group structure of $G$. 

For each $n\in \N$, 
\[ {\cal R}_n \ = \ \{ \ (x,y)\in X\times X:
\mbox{$x$ and $y$ agree except for the initial $n$ edges} \}\]
is a subrelation of ${\cal R}$. Give ${\cal R}_n$ the relative 
topology ${\cal T}_n$ of $X\times X$. 
Then $({\cal R}_n,{\cal T}_n)$ is a CEER and we have 
$({\cal R},{\cal T})=\lim ({\cal R}_n,{\cal T}_n)$. 
Take $\gamma \in G$. As in the proof of the above proposition, 
we can find $n\in \N$ such that 
$(x,\gamma (x))$ is in ${\cal R}_n$ for all $x\in X$. 
The infinite paths $x$ and $\gamma (x)$ agree from the $(n+1)$-st 
level on. Since $\gamma$ is continuous, for sufficiently large 
$m>n$, we may assume that the initial $n$ edges of $\gamma (x)$ 
depend only on the initial $m$ edges of $x$. 
Then, if the initial $m$ edges of $x$ and $y$ coincide, then 
the initial $m$ edges of $\gamma (x)$ and $\gamma (y)$ also coincide, 
because $\gamma$ preserves the $k$-th edges if $k>n$. 
Hence $\gamma$ induces a permutation on the set of paths from 
the top vertex $v_0$ to each $v\in V_m$. These permutations 
form a subgroup $G_m$ in $G$, which is isomorphic to 
$\bigoplus _{v\in V_m}S_{h(v)}$, where $h(v)$ denotes 
the number of paths from the top vertex $v_0$ to $v$. 
Clearly $G_m$ is contained in $G_{m+1}$ and $G$ equals 
the union of all the $G_m$'s. The inclusion $G_m\subset G_{m+1}$ is 
described by the edge set $E_{m+1}$: 
if there are $k$ edges between $v\in V_m$ 
and $v'\in V_{m+1}$, then $S_{h(v')}\subset G_{m+1}$ contains 
$k$ copies of $S_{h(v)}$. 

We can summarize the conclusion just obtained as follows. 
\begin{prop}
When $B=(V,E)$ is a Bratteli diagram and 
$({\cal R},{\cal T})$ is the associated AF equivalence relation, 
the AF full group $G=[{\cal R},{\cal T}]$ can be written 
as a increasing union of subgroups $G_m$, which is isomorphic to 
$\bigoplus _{v\in V_m}S_{h(v)}$. 
\end{prop}
Let $p$ be a path from $v_0$ to $v\in V_m$ and let $U(p)$ be 
a clopen set of infinite paths whose initial paths agree with $p$. 
The zero-dimensional topology in $X$ is generated by clopen sets 
of this form. If $\gamma$ is in $G_m$, then $1_{U(p)}\gamma ^{-1}=
1_{\gamma (U(p))}$ is equal to $U(q)$ for another path $q$ from 
$v_0$ to $v$. Define the $K^0$-group associated with the AF equivalence 
relation $({\cal R},{\cal T})$ by 
\[ K^0({\cal R},{\cal T}) \ = \ C(X, \Z)/\{f-f\gamma ^{-1}:
\gamma \in [{\cal R},{\cal T}] \}. \]
We denote the equivalence class of $f\in C(X,\Z)$ by $[f]$. 
Obviously $K^0({\cal R},{\cal T})$ has a natural unital ordered 
group structure and, by the above observation, it is unital order 
isomorphic to $K_0(B)$. 

Let us consider normal subgroups of $G=[{\cal R},{\cal T}]$. 
The symmetric group $S_{h(v)}$ has the normal subgroup $A_{h(v)}$, 
and so $G_m$ has the normal subgroup $H_m$ isomorphic to 
$\bigoplus _{v\in V_m}A_{h(v)}$. 
\begin{lem}\label{DAFsimple}
In the above setting, the commutator subgroup $D(G)$ equals to 
the union $\bigcup H_m$. Furthermore, $D(G)$ is simple 
if the Bratteli diagram $B$ is simple, 
or equivalently ${\cal R}$ is minimal. 
\end{lem}
\begin{proof}
The first half of the statement is obvious. 

To prove the latter half, let $N$ be a normal subgroup of $D(G)$ 
and let $\gamma\in N$ be a non-trivial element. 
Suppose that $\gamma$ is contained in $H_m$ and its $v$-summand is 
non-trivial for some $v\in V_m$. If $l>m$ is sufficiently large, 
then, for all $w\in V_l$, there exists a path from $v$ to $w$. 
Therefore the $w$-summand of $\gamma$ is non-trivial for all 
$w\in V_l$. 
By taking commutators, we can see that 
every $A_{h(w)}$ is contained in $N$, confirming that 
$H_l$ is contained in $N$, if $l$ is sufficiently large. 
It follows that we have $N=D(G)$. 

Conversely, suppose that $B$ is not simple. We can find a vertex 
$v\in V_m$ such that $V'_l=\{w\in V_l:w\mbox{ is 
connected to }v\}$ is a proper subset of $V_l$ for all $l>m$. 
Let $N_l$ be the direct sum of all $w$-summands for $w\in V'_l$ 
in $H_l$. Then $N=\bigcup_{l>m}N_l$ is a non-trivial normal 
subgroup of $H$. This completes the proof. 
\end{proof}
The quotient group $G_m/H_m$ is isomorphic to $\Z_2^{V_m}$ 
unless $h(v)$ is one for some $v\in V_m$, 
where $h(v)$ denotes the number of paths from the top vertex $v_0$ to $v$.
Suppose that $\#[x]_{\cal R}$ is more than one for all $x\in X$. 
Then $h(v)$ is more than one for all $v\in V_m$, 
if $m$ is sufficiently large. 
Therefore $G/D(G)$ is isomorphic to the inductive limit of 
$G_m/H_m\cong \Z_2^{V_m}$, where the connecting homomorphism 
from $\Z_2^{V_m}$ to $\Z_2^{V_{m+1}}$ is given by the edge set 
$E_m$. Hence we obtain the following 
under the identification of $K^0({\cal R},{\cal T})\otimes \Z_2$ 
with $K^0({\cal R},{\cal T})/2K^0({\cal R},{\cal T})$. 
\begin{lem}\label{AFsignature}
If $\#[x]_{\cal R}$ is more than one for all $x\in X$, 
the quotient group $G/D(G)$ is isomorphic to 
$K^0({\cal R},{\cal T})/2K^0({\cal R},{\cal T})$. 
\end{lem}
\begin{df}\label{sgnonAF}
We denote the composition of the quotient map from $G$ to $G/D(G)$ and 
the isomorphism in the lemma above 
by $\sgn:G\to K^0({\cal R},{\cal T})/2K^0({\cal R},{\cal T})$ and 
call it the signature map. 
For $\gamma \in G=[{\cal R},{\cal T}]$, 
we call $\sgn (\gamma )\in K^0({\cal R},{\cal T})/2K^0({\cal R},{\cal T})$ 
the signature of $\gamma$. 
\end{df}
We would like to give an explicit procedure 
for computing the signature. 
Let $\gamma$ be an element of $[{\cal R},{\cal T}]$. 
For $x\in X$, let $f(x)$ be the least natural number such that 
$\gamma ^{f(x)}(x)=x$. Then $f$ is a well-defined continuous 
function. If $f^{-1}(n)$ is not empty, 
we can find a clopen subset $V_n\subset f^{-1}(n)$ so that 
$V_n, \gamma (V_n), \dots ,\gamma ^{n-1}(V_n)$ are mutually disjoint 
and $f^{-1}(n)=V_n\cup \gamma (V_n)\cup \dots \cup 
\gamma ^{n-1}(V_n)$. Put 
\[ g(x) \ = \ \sum _{n\in 2\N}1_{V_n}(x) \]
for $x\in X$. 
The following lemma claims that $\sgn(\gamma)$ depends 
only on the dynamical behavior of $\gamma$. 
\begin{lem}\label{computesgn}
In the above setting, $\sgn (\gamma)$ is equal to 
$[g]+2K^0({\cal R},{\cal T})$. 
\end{lem}
\begin{proof}
Put $U=\bigcup_{n\in 2\N}V_n$. 
For $z\in X$, we define 
\[ \gamma'(z)=\begin{cases}\gamma(z) & z\in U \\
\gamma^{-1}(z) & z\in\gamma(U) \\
z & \text{otherwise}. \end{cases} \]
Then it is not hard to see that $\gamma'$ belongs to $G$. 
In addition, $\gamma^{-1}\gamma'$ is in $D(G)$, 
and so it suffices to show 
$\sgn(\gamma')=[1_U]+2K^0({\cal R},{\cal T})$. 
Let us assume that 
$\gamma'$ belongs to $G_m\cong\bigoplus_{v\in V_m}S_{h(v)}$. 
For $v\in V_m$, let $P_v$ be the set of paths $p$ 
from $v_0$ to $v$ such that $U(p)\subset U$. 
The $S_{h(v)}$ component of $\gamma'\in G_m$ is 
a product of $\#P_v$ transpositions. 
Hence 
\[ \sgn(\gamma')=\sum_{\#P_v\notin2\N}[1_{U(p_v)}]
+2K^0({\cal R},{\cal T}), \]
where $p_v$ is a path from $v_0$ to $v$. 
It follows that 
\begin{align*}
\sgn(\gamma')&=\sum_{v\in V_m}\#P_v[1_{U(p_v)}]
+2K^0({\cal R},{\cal T}) \\
&=\sum_{v\in V_m}\sum_{p\in P_v}[1_{U(p)}]
+2K^0({\cal R},{\cal T}) \\
&=[1_U]+2K^0({\cal R},{\cal T}), 
\end{align*}
which completes the proof. 
\end{proof}
Note that the construction of the function $g$ is possible, 
if $\gamma$ is of finite order. 

\section{Signatures}
Let $\xp$ be a Cantor minimal system. In this section we would like 
to investigate $[[\phi]]_0$ by using the signature map. 
\begin{lem}\label{decomp}
When $x$ and $y$ have distinct orbits in $X$, we have $[[\phi]]_0=
[[\phi]]_x[[\phi]]_y$, that is, 
for every $\gamma \in [[\phi]]_0$ there exist $\gamma _1 \in 
[[\phi]]_x$ and $\gamma _2\in [[\phi]]_y$ such that 
$\gamma =\gamma _1\gamma _2$. 
\end{lem}
\begin{proof}
Let us denote the forward and backward orbit of $x$ by 
$O^+(x)=\{\phi^n(x);n\in \N\}$ and 
$O^-(x)=\{\phi^{1-n}(x);n\in \N\}$, respectively. 
Put 
\[ A \ = \ \{n\in \N:\gamma (\phi ^n(x))\in O^-(x) \}\]
and
\[ B \ = \ \{n\in \N:\gamma (\phi ^{1-n}(x))\in O^+(x) \},\]
then $A$ and $B$ have the same cardinalities 
as $\gamma \in [[\phi]]_0$, according to \cite[Remark 5.6]{GPS2}. 
Let $l$ be the maximum of $n\in A\cup B$. 
Choose a small clopen neighborhood $U$ of $x$ so that 
$\phi^{-l}(U),\phi ^{-l+1}(U),\dots ,\phi^{l}(U)$ are mutually 
disjoint and 
\[ V \ = \ \bigcup _{n\in B}\phi^{1-n}(U)\cup 
\bigcup _{n\in A}\phi ^n(U) \]
does not contain $y,\phi(y),\dots ,\phi ^{2l-1}(y)$. 
Let $\pi$ be a bijective map from $A$ to $B$. 
Define $\gamma _2\in \Homeo(X)$ by 
\[ \gamma _2(z) \ = \ \begin{cases}
\phi^{1-\pi(n)-n}(z) & z\in \phi^n(U)\mbox{ for some }n\in A \\
\phi^{n-1+\pi^{-1}(n)}(z) & z\in \phi ^{1-n}(U)
\mbox{ for some }n\in B \\
z & z\in X\setminus V \end{cases}. \]
Then, it is not hard to see 
$\gamma _2\in [[\phi]]_y$ and $\gamma\gamma_2\in [[\phi]]_x$. 
\end{proof}
The signature map has been already defined on $[[\phi]]_x$ 
for every $x\in X$, 
taking its value in $K^0({\cal R}_x, {\cal T}_x)\otimes \Z_2$. 
Since $K^0({\cal R}_x, {\cal T}_x)$ can be identified with 
$K^0\xp$ (see \cite[Theorem 3.7]{GPS1}), we may regard 
$\sgn(\gamma)\in K^0\xp\otimes \Z_2$. 

We need to know the behavior of the signature map 
on the intersection $[[\phi ]]_x\cap [[\phi ]]_y$. 
\begin{lem}\label{intersect1}
For $\gamma \in [[\phi]]_x\cap [[\phi]]_y$, the signature 
$\sgn (\gamma )\in K^0\xp/2K^0\xp$ is well-defined. 
\end{lem}
\begin{proof}
It is clear from Lemma \ref{computesgn}. 
\end{proof}
\begin{lem}\label{intersect2}
The signature map $\sgn$ from $[[\phi]]_x\cap [[\phi]]_y$ to 
$K^0\xp/2K^0\xp$ is surjective. 
\end{lem}
\begin{proof}
The subequivalence relation ${\cal R}_{\{x,y\}}$ is also an AF 
equivalence relation by Theorem \ref{AFinCM} and 
$[[\phi]]_x\cap [[\phi]]_y=
[{\cal R}_{\{x,y\}},{\cal T}_{\{x,y\}}]$. 
The image of the signature map defined on 
$[{\cal R}_{\{x,y\}},{\cal T}_{\{x,y\}}]$ is 
the whole of $K^0({\cal R}_{\{x,y\}},{\cal T}_{\{x,y\}})
\otimes \Z_2$ by Lemma \ref{AFsignature}. 
There exists a canonical surjective map from 
$K^0({\cal R}_{\{x,y\}},{\cal T}_{\{x,y\}})$ to $K^0\xp$ 
(see \cite[Theorem 1.17]{GPS1}), 
thereby completing the proof. 
\end{proof}
\begin{lem}
Let $x$ and $y$ be lying in distinct orbits. 
For $\gamma \in [[\phi]]_0$, let $\gamma =\gamma _1\gamma _2$ be 
the decomposition described in Lemma \ref{decomp}, 
where $\gamma _1\in [[\phi]]_x$ and $\gamma _2\in [[\phi]]_y$. 
Put $s=\sgn(\gamma_1)+\sgn(\gamma_2)\in K^0\xp/2K^0\xp$. 
\begin{enumerate}
\item $s$ does not depend on the choice of $\gamma_1$ and $\gamma_2$. 
\item $s$ does not depend on the choice of $x$ and $y$ 
lying in distinct orbits. 
\end{enumerate}
\end{lem}
\begin{proof}
(1) follows from Lemma \ref{intersect1} and 
the fact that $\sgn$ is a homomorphism on AF full groups. 
Let us prove (2). 
Suppose that $x,y$ and $z$ have distinct orbits. In the proof of Lemma 
\ref{decomp} it is easy to construct $\gamma_2$ so that 
$\gamma_2$ does not move $z,\phi(z),\dots ,\phi ^{l-1}(z)$. 
Therefore we can make $\gamma_2$ in $[[\phi]]_y\cap [[\phi]]_z$, 
which means that the signature takes the same value if one uses 
$x$ and $z$ instead of $x$ and $y$ for the definition. 
In the case of using $z$ and $y$ to define the signature, 
we can show the assertion similarly. 
\end{proof}
\begin{df}
Let $\gamma\in[[\phi]]_0$. 
We denote the value $s\in K^0\xp/2K^0\xp$ 
described in the lemma above by $\sgn(\gamma)$. 
By the lemma above, $\sgn$ is a well-defined map 
from $[[\phi]]_0$ to $K^0\xp/2K^0\xp$. 
We call $\sgn$ the signature map. 
\end{df}
\begin{prop}
The signature map from $[[\phi]]_0$ to $K^0\xp/2K^0\xp$ is 
a group homomorphism. 
\end{prop}
\begin{proof}
It suffices to show $\sgn(\gamma \tau)=\sgn(\gamma)+\sgn(\tau)$ 
for $\gamma ,\tau\in [[\phi]]_0$. Suppose that $x,y\in X$ have 
distinct orbits. Let $\tau=\tau_1\tau_2$ be the decomposition 
as in Lemma \ref{decomp}. Put 
\[ m \ = \ \max\{n\in \N:\tau _1^{-1}(\phi ^{n}(y))\in 
O^-(y)\mbox{ or }\tau _1^{-1}(\phi ^{1-n}(y))\in O^+(y) \}.\]
When we apply Lemma \ref{decomp} to $\gamma$, 
it is possible to choose 
the clopen set $U$ so that $V$ does not contain $\phi ^k(y)$ 
for all $1-m\leq k\leq m$. Then, the obtained 
$\gamma_2\in [[\phi]]_y$ satisfies $\tau_1^{-1}\gamma_2\tau_1
\in [[\phi]]_y$. 
Because the signature map defined on an AF full group is 
a homomorphism (Definition \ref{sgnonAF}), we have 
$\sgn (\gamma_1\tau_1)=\sgn(\gamma_1)+\sgn(\tau_1)$. 
Similarly we get $\sgn((\tau_1^{-1}\gamma_2\tau_1)\tau_2)=
\sgn(\tau_1^{-1}\gamma_2\tau_1)+\sgn(\tau_2)=\sgn(\gamma_2)+
\sgn(\tau_2)$, where the last equation follows from Lemma 
\ref{computesgn}. Hence we have 
\begin{align*}
\sgn(\gamma\tau) &= \sgn(\gamma_1\tau_1)+
\sgn(\tau_1^{-1}\gamma_2\tau_1\tau_2) \\
&= \sgn(\gamma_1)+\sgn(\tau_1)+\sgn(\gamma_2)+\sgn(\tau_2) \\
&= \sgn(\gamma)+\sgn(\tau), 
\end{align*}
proving that $\sgn$ is a homomorphism. 
\end{proof}
\begin{lem}\label{dense}
Let $U$ be a clopen neighborhood of $x\in X$. 
Then there exists a clopen neighborhood $V$ of $x$ such that 
$V\subset U$ and $[1_V]$ is $2$-divisible in $K^0\xp$. 
\end{lem}
\begin{proof}
Choose a small neighborhood $V_0$ contained in $U$ so that 
$2\mu(V_0)$ is less than $\mu(U)$ for all invariant probability 
measure $\mu$. 
By applying \cite[Lemma 2.5]{GW} with $A=U\setminus V_0$ and 
$B=V_0$, we obtain $V_1\subset U\setminus V_0$ satisfying 
$[1_{V_0}]=[1_{V_1}]$ in $K^0\xp$. 
Then $V_0\cup V_1$ does the work. 
\end{proof}
Now we are ready to prove the main theorems of this section. 
\begin{thm}\label{CMsignature}
Let $\xp$ be a Cantor minimal system. Then the kernel of 
the signature map equals the commutator $D([[\phi]]_0)$. 
Thus $\sgn$ establishes an isomorphism between 
$[[\phi]]_0/D([[\phi]]_0)$ and $K^0\xp/2K^0\xp$. 
\end{thm}
\begin{proof}
It is clear that the commutator subgroup $D([[\phi]]_0)$ is 
contained in the kernel of $\sgn$. Let us prove the other 
inclusion. Take $\gamma \in [[\phi]]_0$ with $\sgn(\gamma)=0$. 
In the proof of Lemma \ref{decomp} it is easy to choose the clopen 
set $U$ so that $[1_U]$ is 2-divisible in $K^0\xp$ 
by using Lemma \ref{dense}. 
Then, by Lemma \ref{computesgn}, we can see that $\sgn(\gamma_2)=0$. 
As $\sgn(\gamma)=0$, the signature of $\gamma_1$ is equal to zero, 
too. Lemma \ref{AFsignature} tells us that $\gamma_1$ and $\gamma_2$ 
are contained in the commutator subgroups of AF full groups. 
Hence $\gamma$ is in the commutator subgroup $D([[\phi]]_0)$. 

Since the signature map is surjective, it is easily verified that 
$[[\phi]]_0/D([[\phi]]_0)$ is isomorphic to $K^0\xp/2K^0\xp$. 
\end{proof}
Note that $\sgn(\phi\gamma\phi^{-1})$ is equal to $\sgn(\gamma)$ 
for all $\gamma \in D([[\phi]]_0)$. 
It follows that the commutator subgroup $D([[\phi]]_0)$ 
coincides with $D([[\phi]])$. 
\begin{thm}\label{CMsimple}
Let $\xp$ be a Cantor minimal system. 
The commutator subgroup $D([[\phi]]_0)$ is simple. 
\end{thm}
\begin{proof}
Suppose that $N\subset D([[\phi]]_0)$ is a nontrivial normal 
subgroup. 
It suffices to show that $N$ contains $D([[\phi]]_x)$ for every 
$x\in X$. 
Take a nontrivial element $\gamma \in N$. 
Let $l$ be the maximum of $n(z)$, where $n$ is the continuous 
function from $X$ to $\Z$ satisfying $\gamma(z)=\phi^{n(z)}(z)$. 
In view of $\gamma \neq \id ,\phi$ we can select $y\in X$ 
so that $x$ and $y$ have distinct orbits and 
$\gamma(y)\neq y,\phi(y)$. 
Choose a clopen neighborhood $U$ of $y$ 
which does not intersect with $\phi(U)$ and does not 
contain $\phi(y),\gamma(y),\phi^{-1}\gamma(y)$ and 
$\phi^{k}(x)$ with $|k|\leq l+1$. 
From Lemma \ref{dense} 
we may further assume that $[1_U]$ is 2-divisible in $K^0\xp$. 
Define $\sigma\in D([[\phi]]_0)$ by 
\[ \sigma(z) \ = \ \begin{cases}
\phi(z) & z\in U \\
\phi^{-1}(z) & z\in \phi(U) \\
z & \mbox{otherwise}. \end{cases}\]
From the construction, we infer that 
the commutator $\sigma\gamma^{-1}\sigma\gamma$ belongs to 
$\ker\sgn\cap[[\phi]]_x=D([[\phi]]_x)$. 
Besides, on account of $\sigma\gamma^{-1}\sigma\gamma(y)=
\sigma\gamma^{-1}\gamma(y)=\sigma(y)=\phi(y)$, 
the commutator $\sigma\gamma^{-1}\sigma\gamma$ is not the identity. 
Then using Lemma \ref{DAFsimple} shows that 
$N$ contains $D([[\phi]]_x)$. 
\end{proof}
Combining Theorem \ref{CMsignature} and \ref{CMsimple}, 
we get the following. 
\begin{cor}\label{CMsimpleCOR}
For a Cantor minimal system $\xp$, 
$[[\phi]]_0$ is simple if and only if 
$K^0\xp$ is $2$-divisible. 
\end{cor}

We would like to conclude this section by giving an example 
of a computation of the signature. 
Let $U$ be a nonempty clopen subset of $X$ and 
let $\phi_U$ be the first return map on $U$. 
We can extend $\phi_U$ to a homeomorphism on $X$ 
by putting $\phi_U(x)=x$ for all $x\in X\setminus U$. 
We remark that the index of $\phi_U$ is one. 
\begin{lem}
When $U$ and $V$ are nonempty clopen subsets, the signature of 
$\phi_U\phi_V^{-1}$ is equal to $[1_U]+[1_V]+2K^0\xp$. 
\end{lem}
\begin{proof}
We may assume that $V$ is contained in $U$, 
because $\phi_U\phi_V^{-1}=
(\phi_{U\cup V}\phi_U^{-1})^{-1}(\phi_{U\cup V}\phi_V^{-1})$. 
Then, it is easily checked that $\gamma=\phi_U\phi_V^{-1}$ is 
in $[[\phi]]_x$ for every $x\in V$. 
For $x\in X$, let $f(x)$ be the least natural number 
such that $\gamma^{f(x)}(x)=x$. Put $V_n=V\cap f^{-1}(n)$. 
Since 
\[ U=
\bigcup_{n\in \N}\bigcup_{k=1}^n\gamma^k(V_n), \]
by virtue of Lemma \ref{computesgn}, we obtain 
\begin{align*}
\sgn(\gamma) &= \sum _{n\in 2\N}[1_{V_n}]+2K^0\xp \\
&= \sum _{n\in 2\N}[1_{V_n}]+
\sum_{n\in 2\N}\sum_{k=1}^n[1_{\gamma^k(V_n)}]+
\sum_{n\in 2\N-1}\left([1_{V_n}]+\sum_{k=1}^n[1_{\gamma^k(V_n)}]
\right)+2K^0\xp \\
&= \sum _{n\in \N}[1_{V_n}]+
\sum_{n\in \N}\sum_{k=1}^n[1_{\gamma^k(V_n)}]+2K^0\xp \\
&= [1_V]+[1_U]+2K^0\xp .
\end{align*}
\end{proof}

\section{Finite generatedness}
Let $\xp$ be a Cantor minimal system. 
In this section we will give a necessary and sufficient condition 
for $D([[\phi]]_0)$ to be finitely generated. 
As pointed out in the proof of Theorem \ref{CMsignature}, 
every $\gamma\in D([[\phi]]_0)$ is a product of some 
$\gamma_1\in D([[\phi]]_x)$ and $\gamma_2\in D([[\phi]]_y)$. 
Thus, in order to prove the finite generatedness, 
it suffices to show that 
there exists a finite number of elements of $D([[\phi]]_0)$ 
generating $D([[\phi]]_x)$ for all $x\in X$. 
\begin{lem}\label{alter}
Let $n$ be a natural number more than two. 
The alternating group $A_n$ is generated by 
$(1,2,3),(2,3,4),\dots ,(n-2,n-1,n)$. 
\end{lem}
\begin{proof}
The proof is by induction. 
For $n=3,4$ the assertion follows from the straightforward 
computation. 
Let us assume that $A_n$ is generated by 
$(1,2,3),(2,3,4),\dots ,(n-2,n-1,n)$. 

Suppose $k\in \{1,2,\dots ,n-2\}$. Then 
\[ ((k,n-1)(n-1,n))(n-1,n,n+1)((k,n-1)(n-1,n))^{-1}
\ = \ (n,k,n+1). \]
Similarly $(j,k,n+1)$ can be generated for all distinct 
$j,k\in \{1,2,\dots ,n\}$. 
\end{proof}

Let $U$ be a nonempty clopen subset of $X$ such that 
$\phi ^{-1}(U),U$ and $\phi (U)$ are mutually disjoint. 
Then we define a homeomorphism $\gamma _U$ of order three by 
\[ \gamma _U(x) \ = \ \begin{cases}
\phi (x) & x\in \phi ^{-1}(U)\cup U \\
\phi ^{-2}(x) & x\in \phi (U) \\
x & \mbox{otherwise}. \end{cases} \]
For convenience, we put $\gamma_{\emptyset}=\id$. 

When $\phi ^{-2}(U),\phi ^{-1}(U),U,\phi (U)$ and $\phi ^2(U)$ are 
mutually disjoint, 
we define a homeomorphism $\tau _U$ of order five by 
$\tau _U=\gamma _{\phi^{-1}(U)}\gamma _{\phi (U)}$. 
Observe that all $\gamma _U$ and $\tau _U$ belong to 
the commutator subgroup $D([[\phi]]_0)$, 
because they have odd order (see Lemma \ref{computesgn}). 
\begin{lem}\label{gamU}
The commutator subgroup $D([[\phi]]_0)$ is generated by 
\[ \Gamma=\{ \gamma _U\in D([[\phi]]_0):\phi ^{-1}(U),U\mbox{ and }
\phi(U)\mbox{ are mutually disjoint} \}.\]
\end{lem}
\begin{proof}
It suffices to show that for all $x\in X$, 
$D([[\phi]]_x)$ is generated by elements of $\Gamma$. 
Lemma \ref{DAFsimple} tells us that $D([[\phi]]_x)$ is 
the inductive limit of finite groups which are isomorphic to 
direct sums of alternating groups. 
In these alternating groups, every 3-cyclic permutation of the form 
$(i-1,i,i+1)$ corresponds to some $\gamma_U$. 
Thanks to Lemma \ref{alter}, we get the conclusion. 
\end{proof}

\begin{lem}\label{gamUtauV}
Let $U$ and $V$ be clopen subsets of $X$. 
\begin{enumerate}
\item If $\phi^{-2}(V),\phi ^{-1}(V),V,\phi (V)$ and $\phi ^2(V)$ 
are mutually disjoint and $U$ is contained in $V$, then 
we have $\tau_V\gamma_U\tau_V^{-1}=\gamma_{\phi(U)}$ and 
$\tau_V^{-1}\gamma_U\tau_V=\gamma_{\phi^{-1}(U)}$. 
\item If $\phi^{-1}(U),U,\phi(U)\cup \phi^{-1}(V),V$ and $\phi(V)$ 
are mutually disjoint, then we have 
$\gamma_V\gamma_U^{-1}\gamma_V^{-1}\gamma_U=
\gamma_{\phi(U)\cap \phi^{-1}(V)}$. 
\end{enumerate}
\end{lem}
\begin{proof}
Straightforward computations. 
\end{proof}

The following is the main theorem of this section. 
The reader, who is not familiar with the rudiments of the theory 
of subshifts, may refer to \cite{LM}. 
\begin{thm}\label{Dfg}
Let $\xp$ be a Cantor minimal system. 
The commutator subgroup $D([[\phi]]_0)$ is finitely generated 
if and only if $\xp$ is conjugate to a minimal subshift. 
\end{thm}
\begin{proof}
The `only if' part is easy to show. Let us assume that 
$D([[\phi]]_0)$ is generated by finitely many elements 
$\gamma_1,\gamma_2,\dots ,\gamma_m$. 
Let $n_i:X\rightarrow \Z$ be the continuous function such that 
$\gamma_i(x)=\phi^{n_i(x)}(x)$ for all $x\in X$. 
Then, ${\cal P}_i=\{n_i^{-1}(k)\}_{k\in \Z}$ is 
a finite clopen partition for each $i=1,2,\dots ,m$. 
Take a finite clopen partition ${\cal P}$ which is finer than 
all the ${\cal P}_i$'s. 
We equip ${\cal P}^{\Z}$ with the product topology and 
denote the shift map by $\sigma$. 
Define a continuous map $\pi$ from $X$ to ${\cal P}^{\Z}$ 
so that $\phi^k(x)\in \pi(x)_k$. 
Then $\pi$ is a factor map from $\xp$ to the minimal subshift 
$(\pi(X),\sigma)$. 
For each $i=1,2,\dots ,m$, we can define a homeomorphism 
$\tau_i\in \Homeo(\pi(X))$ by $\tau_i(\xi) \ = \ \sigma^k(\xi)$ 
when $\xi _0\subset n_i^{-1}(k)$. 
It is an easy observation that $\tau_i$ is an element of $[[\sigma]]$ 
and $\pi\gamma_i=\tau_i\pi$. 
We would like to prove that $\pi$ is injective. 
Suppose that $x$ and $y$ are distinct points in $X$ and 
$\pi(x)=\pi(y)$. There exists $\gamma\in D([[\phi]]_0)$ such that 
$\gamma(x)\neq x$ and $\gamma(y)=y$. Since $D([[\phi]]_0)$ is 
generated by $\gamma_1,\gamma_2,\dots ,\gamma_m$, 
we can find $r_1,r_2,\dots ,r_l\in \Z$ and 
$i_1,i_2,\dots ,i_l\in \{1,2,\dots ,m\}$ such that 
$\gamma =\gamma_{i_1}^{r_1}\gamma_{i_2}^{r_2}\dots\gamma_{i_l}^{r_l}$. 
This together with $\pi\gamma_i=\tau_i\pi$ shows that 
\begin{align*}
\pi\gamma(x) &= 
\pi\gamma_{i_1}^{r_1}\gamma_{i_2}^{r_2}\dots\gamma_{i_l}^{r_l}(x) \ 
= \ \tau_{i_1}^{r_1}\tau_{i_2}^{r_2}\dots\tau_{i_l}^{r_l}\pi(x) \\
&= \tau_{i_1}^{r_1}\tau_{i_2}^{r_2}\dots\tau_{i_l}^{r_l}\pi(y) \ 
= \ \pi\gamma_{i_1}^{r_1}\gamma_{i_2}^{r_2}\dots\gamma_{i_l}^{r_l}(y) \\
&= \pi\gamma (y) \ = \ \pi(y) \ = \ \pi(x), 
\end{align*}
which means $\sigma^k\pi(x)=\pi\phi^k(x)=\pi(x)$ for some nonzero 
$k\in \Z$. This is a contradiction because $\sigma$ is minimal. 

Let us prove the `if' part. Suppose that $\xp$ is conjugate to 
a minimal subshift. We may assume that $X$ is a shift invariant closed 
subset of $A^\Z$, where $A$ is a finite set. 
We may further assume that $\xi_i\neq \xi_j$ for all $\xi\in X$ and 
$i,j\in \Z$ with $|i-j|\leq 4$. 
We define the cylinder sets by 
\[ [a_{-m}\dots a_{-1}a_0.a_1\dots a_n] \ 
= \ \{\xi \in X:\xi_i=a_i\mbox{ for all }
i=-m,\dots ,n-1,n\} \]
for $n,m\in \N\cup\{0\}$ and $a_{-m},\dots ,a_n\in A$. 
Note that $\phi^{-2}(U),\phi^{-1}(U),U,\phi(U)$ and $\phi^2(U)$ are 
mutually disjoint, if $U$ is a cylinder set. 
Since every clopen set can be written 
as a disjoint union of these cylinder sets, 
it suffices to show that $\gamma_U$ can be generated for every 
cylinder set $U$ by virtue of Lemma \ref{gamU}. 
Put 
\[ F \ = \ \{\gamma_U:U=[ab.c], \ a,b,c\in A\}. \]
We claim that $F$ generates $D([[\phi]]_0)$. 
It is not hard to see that $\tau_{[a.]}$ can be generated 
for all $a\in A$, 
because $\gamma_{\phi^{-1}([a.])}$ and $\gamma_{\phi([a.])}$ 
are generated. 
When $U$ is a cylinder set $[a_{-m}\dots a_0.\dots a_n]$, 
by Lemma \ref{gamUtauV} (i), it is easily verified that 
\[ \tau_{[a_0.]}\gamma_U\tau_{[a_0.]}^{-1} \ = 
\ \gamma_{\phi(U)} \]
and 
\[ \tau_{[a_0.]}^{-1}\gamma_U\tau_{[a_0.]} \ = 
\ \gamma_{\phi^{-1}(U)}. \]
Therefore, it suffices to show that $\gamma_U$ can be generated 
for every cylinder set $[a_{-m}\dots a_0.a_1]$. 
This is inductively done by using Lemma \ref{gamUtauV} (ii) with 
$U=[a_{-m}\dots a_0.a_1]$ and $V=[a_1a_2.]$. 
\end{proof}
This result allows us to conclude the following. 
\begin{cor}\label{DfgCOR}
When $\xp$ is a Cantor minimal system, the following are equivalent. 
\begin{enumerate}
\item $\xp$ is a minimal subshift and $K^0\xp \otimes \Z_2$ is 
a finite group. 
\item $[[\phi]]_0$ is finitely generated. 
\item $[[\phi]]$ is finitely generated. 
\end{enumerate}
\end{cor}
\begin{proof}
(i)$\Rightarrow$(ii). It follows from Theorem \ref{CMsignature} and 
Theorem \ref{Dfg}. 

(ii)$\Rightarrow$(iii). It follows from 
$[[\phi]]/[[\phi]]_0\cong \Z$. 

(iii)$\Rightarrow$(i). It is clear that 
$[[\phi]]/D([[\phi]]_0)\cong \Z\oplus (K^0\xp\otimes \Z_2)$ 
should be finitely generated. 
By the same proof as the `only if' part of Theorem \ref{Dfg}, 
we can conclude that $\xp$ is a minimal subshift. 
\end{proof}
There exist uncountably many minimal subshifts $\xp$ 
whose $K^0$-group $K^0\xp$ is $2$-divisible by \cite[Theorem 1.1]{Su}. 
Combining Theorem \ref{fulliso} (iii), Corollary \ref{CMsimpleCOR} 
and the result obtained above, 
we get uncountably many finitely generated simple groups. 
\bigskip

We would like to show that 
the topological full group is never finitely presented. 
In the following lemma, topological entropy is 
denoted by $h(\cdot)$. 
A dynamical system $\xp$ is said to be embeddable into 
another dynamical system $\yp$, 
if there exists a continuous injection $\pi:X\rightarrow Y$ 
with $\pi\phi=\psi\pi$. 
\begin{lem}
Let $\xp$ be a minimal subshift. 
Let $\yp$ be an irreducible subshift of finite type 
and let $p\in \N$ be its period. 
Then $\xp$ can be embedded into $\yp$ if and only if 
$h(\phi)<h(\psi)$ and $[1_X]$ is divisible by $p$ in $K^0\xp$. 
\end{lem}
\begin{proof}
When $p$ is one, this is a direct consequence of 
Corollary 10.1.9 of \cite{LM}. 
Let $p\geq 2$. 
There exists a clopen $\psi^p$-invariant subset $U$ of $Y$ and 
$Y$ is the disjoint union of $U,\psi(U),\dots,\psi^{p-1}(U)$. 

Suppose that $\xp$ is embedded into $\yp$. 
Evidently we have $h(\phi)<h(\psi)$ 
(see \cite[Corollary 4.4.9]{LM}). 
Define $V=U\cap X$. Then $[1_X]$ is equal to $p[1_V]$. 

Let us prove the converse. 
We can find a clopen subset $V\subset X$ such that 
$X$ is the disjoint union of $V,\phi(V),\dots,\phi^{p-1}(V)$, 
because $[1_X]$ is divisible by $p$ in $K^0\xp$. 
Then $(V,\phi^p)$ is a minimal subshift. 
Since $(U,\psi^p)$ is a mixing shift of finite type and 
$h(\phi^p)=ph(\phi)<ph(\psi)=h(\psi^p)$, by Corollary 10.1.9 of 
\cite{LM}, $(V,\phi^p)$ is embedded into $(U,\psi^p)$. 
Therefore $\xp$ is embedded into $\yp$. 
\end{proof}
In the proof of the following theorem, we use the fact that 
the commutator subgroup of the topological full group is 
a complete invariant for flip conjugacy. 
This fact was mentioned in \cite[Theorem 4.11]{S2} without proof. 
\begin{thm}
Let $\xp$ be a minimal subshift. 
Then $D([[\phi]]_0)$ is never finitely presented. 
\end{thm}
\begin{proof}
Suppose that $X$ is a subshift in $A^\Z$, where $A$ is a finite set.  
We will freely use the notation of Theorem \ref{Dfg}. 
We may assume that $X$ is contained in 
\[ Y \ = \ \{\xi\in A^\Z:\xi_i\neq\xi_j\mbox{ for all }
i,j\in\Z\mbox{ with }|i-j|\leq4\}. \]
Notice that $Y$ is a subshift of finite type. 
To avoid the ambiguity, denote the shift map on $Y$ by $\phi$ and 
its restriction to a shift invariant subspace $Z$ by $\phi|_Z$. 
We define the cylinder sets as clopen subsets of $Y$. 
If $U$ is a clopen set of $Y$ and $U$, $\phi(U)$, $\phi^2(U)$ 
are disjoint, then we define $\gamma_U\in\Homeo(Y)$ in the same way 
as described before. 
Put 
\[ F \ = \ \{\gamma_U\in\Homeo(Y):
U=[ab.c], \ a,b,c\in A\}. \]

The proof is by contradiction. Suppose that $D([[\phi|_X]]_0)$ is 
finitely presented. 
According to the proof of Theorem \ref{Dfg}, $D([[\phi|_X]]_0)$ 
is generated by $\{\gamma|_X:\gamma\in F\}$, 
and moreover, if $Z$ is another minimal subshift in $Y$, then 
$D([[\phi|_Z]]_0)$ is also generated 
by $\{\gamma|_Z:\gamma\in F\}$. 
By \cite[V.2]{H}, $D([[\phi|_X]]_0)$ may be assumed to have 
a finite presentation 
$\langle(\gamma|_X)_{\gamma\in F}, (r_i)_{i\in I}\rangle$, 
where $I$ is a finite set and 
$r_i$ is an element of the free group over $F$. 
Let $\tilde{r_i}$ denote the canonical image of $r_i$ 
in $\Homeo(Y)$. It is not hard to see that the fixed point set 
$\{y\in Y:\tilde{r_i}(y)=y\}$ is clopen in $Y$ 
for all $i\in I$. 
Since $\tilde{r_i}|_X=\id$, the clopen subset 
\[ Y_0 \ = \ \{y\in Y:
\tilde{r_i}(y)=y\mbox{ for all }i\in I\} \]
contains $X$. 
Hence $X$ is a shift invariant closed subset of 
$Y_1=\bigcap_{n\in\Z}\phi^n(Y_0)$, 
which is also a shift of finite type. 
By taking the irreducible component containing $X$, 
we may further assume that $Y_1$ is irreducible. 
Let $Z$ be any minimal subshift in $Y_1$. 
Then $D([[\phi|_Z]]_0)$ is a quotient of $D([[\phi|_X]]_0)=
\langle(\gamma|_X)_{\gamma\in F}, (r_i)_{i\in I}\rangle$, 
because $\tilde{r_i}|_Z=\id$. 
Therefore we have $D([[\phi|_Z]]_0)\cong D([[\phi|_X]]_0)$ 
by Theorem \ref{CMsimple}, and so 
$\phi|_Z$ is flip conjugate to $\phi|_X$. 
This contradicts the lemma above. 
\end{proof}
\begin{cor}
Let $\xp$ be a minimal subshift. 
Then $[[\phi]]$ and $[[\phi]]_0$ are never finitely presented. 
\end{cor}
\begin{proof}
If $K^0\xp/2K^0\xp$ is infinite, we have nothing to do 
because of Corollary \ref{DfgCOR}. 
Suppose that $K^0\xp/2K^0\xp$ is finite. 
Combining the theorem above and V.5.(i) of \cite{H}, 
we see that $[[\phi]]_0$ is not finitely presented. 
As for $[[\phi]]$, we get the conclusion from \cite[V.15]{H}. 
\end{proof}

\section{Examples}
\begin{exm}
We will consider an example of Cantor minimal systems 
represented by a stationary ordered Bratteli diagram, 
in other words, minimal subshifts arising from substitution rules. 
We refer the reader to \cite{DHS} for the details of these topics. 

Let $B=(V,E)$ be the stationary Bratteli diagram 
given by the following figure, that is, every vertex set 
$V_n$ except for $V_0$ consists of two vertices, 
every edge set $E_n$ except for $E_1$ is 
a copy of the finite set $\{a,b,c,\dots ,h\}$, and 
the range and source map are defined in the same way 
at every level. 

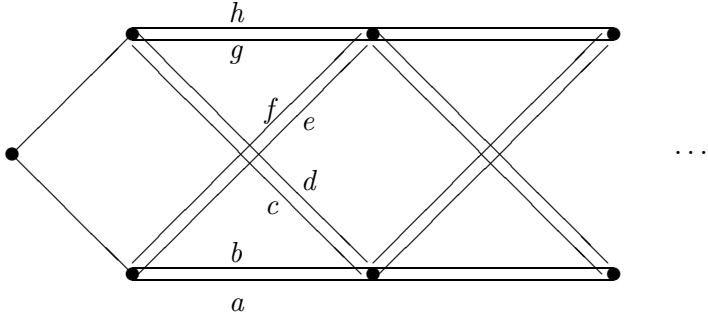
\setlength{\unitlength}{0.8mm}
\begin{picture}(200,60)

\put(10,30){\circle*{2}}
\multiput(30,50)(40,0){3}{\circle*{2}}
\multiput(30,10)(40,0){3}{\circle*{2}}

\put(10,30){\line(1,1){20}}
\put(10,30){\line(1,-1){20}}

\multiput(30,51)(40,0){2}{\line(1,0){40}}
\multiput(30,49)(40,0){2}{\line(1,0){40}}
\multiput(30,11)(40,0){2}{\line(1,0){40}}
\multiput(30,9)(40,0){2}{\line(1,0){40}}

\multiput(31,50)(40,0){2}{\line(1,-1){38}}
\multiput(30,48)(40,0){2}{\line(1,-1){38}}
\multiput(30,12)(40,0){2}{\line(1,1){38}}
\multiput(31,10)(40,0){2}{\line(1,1){38}}

\put(46,4){\textit{a}}
\put(46,12){\textit{b}}
\put(52,20){\textit{c}}
\put(58,24){\textit{d}}
\put(58,34){\textit{e}}
\put(52,36){\textit{f}}
\put(46,46){\textit{g}}
\put(46,52){\textit{h}}

\put(120,30){\ldots}

\end{picture}

We give an ordering on $E_n$ for $n\geq 2$ by 
$a<b<c<d$ and $e<g<f<h$. 
Then $B=(V,E)$ is a stationary ordered Bratteli diagram 
with the unique minimal path $(a,a,\dots )$ and 
the unique maximal path $(h,h,\dots )$. 
Let $\xp$ be the associated Cantor minimal system. 
This Cantor minimal system $\xp$ was presented 
in the final section of \cite{S1} by Skau in order to show that 
there are two topologically orbit equivalent dynamical systems 
that are not conjugate. 
Actually, $\xp$ is strong orbit equivalent to the odometer system 
of type $2^\infty$, because $K^0\xp$ is isomorphic to $\Z[1/2]$ 
(see \cite[Theorem 2.1]{GPS1}). But $\xp$ is not isomorphic to 
the odometer system. 

Since $K^0\xp\cong \Z[1/2]$ is 2-divisible, 
$[[\phi]]_0$ is a simple group by Corollary \ref{CMsimpleCOR}. 
Moreover, Proposition 16 of \cite{DHS} tells us that 
$\xp$ is isomorphic to the minimal subshift arising from 
the substitution rule $0\mapsto 0011$ and $1\mapsto 0101$, and so 
$[[\phi]]_0$ is a finitely generated group by Corollary \ref{DfgCOR}. 

To specify generators of $[[\phi]]_0$, 
we define a clopen set $U(a)$ in the infinite path space $X$ by 
\[ U(a) \ = \ \{(e_n)_{n\in \N}\in X:e_2=a \} \]
and put 
\[ \sigma_a(x) \ = \ \begin{cases}
\phi(x) & x\in U(a) \\
\phi^{-1}(x) & x\in \phi(U(a)) \\
x & \mbox{otherwise}. \end{cases}\]
Clopen sets $U(b),U(c),\dots ,U(h)$ and 
homeomorphisms $\sigma_b,\sigma_c,\dots ,\sigma_h$ are defined 
in a similar fashion. 
Let $\sigma_1=\sigma_a\sigma_e$, $\sigma_2=\sigma_b\sigma_g$, 
$\sigma_3=\sigma_c\sigma_f$ and $\sigma_4=\sigma_d\sigma_h$. 

We infer that $[[\phi]]_0$ is generated by four elements 
$\sigma_1,\sigma_1\sigma_2\sigma_3, \sigma_4$ and $\sigma_d$. 
In order to check it, let $G$ be the subgroup generated 
by this four elements. By the equations 
\begin{align*}
& (\sigma_1\sigma_2\sigma_3)\sigma_1(\sigma_1\sigma_2\sigma_3)^{-1}
=\sigma_2, &&
(\sigma_1\sigma_2\sigma_3)\sigma_2(\sigma_1\sigma_2\sigma_3)^{-1}
=\sigma_3, \\
& (\sigma_2\sigma_3\sigma_4)^{-1}\sigma_d(\sigma_2\sigma_3\sigma_4)
=\sigma_c, &&
(\sigma_2\sigma_3\sigma_4)^{-1}\sigma_c(\sigma_2\sigma_3\sigma_4)
=\sigma_b, \\
& (\sigma_1\sigma_2\sigma_3)^{-1}\sigma_b(\sigma_2\sigma_3\sigma_4)
=\sigma_a, &&
\sigma_e=\sigma_a\sigma_1, \sigma_g=\sigma_b\sigma_2, 
\sigma_f=\sigma_c\sigma_3, \sigma_h=\sigma_d\sigma_4, 
\end{align*}
we get $\sigma_a,\sigma_b,\dots ,\sigma_h\in G$. 
Besides, by $\sigma_a\sigma_4\sigma_a\sigma_4=\gamma_{U(a)}$, 
we have $\gamma_{U(a)}\in G$. In the same way 
$\gamma_{U(b)},\gamma_{U(c)},\dots ,\gamma_{U(h)}$ are also obtained. 
As shown in Lemma \ref{gamU}, to prove $G=[[\phi]]_0$, 
it suffices to show that every $\gamma_U$ is contained in $G$. 
In order to get $\gamma_U$ for a smaller clopen set $U$, 
we use the similar method to the proof of Theorem \ref{Dfg}. 
For example, $\gamma_{U(a)\cap \phi(U(d))}$ is obtained 
as the commutator $\gamma_{U(b)}\sigma_d\gamma_{U(b)}^{-1}\sigma_d$. 
We leave the details to the reader. 

Evidently, the topological full group $[[\phi]]$ is generated by 
$[[\phi]]_0$ and $\phi$. 
Note that $[[\phi]]$ has three generators 
$\phi, \sigma_1$ and $\sigma_d$, 
because we have $\phi^i\sigma_1\phi^{-i}=\sigma_{i+1}$ for $i=1,2,3$. 
\end{exm}

\begin{exm}
Let $\alpha\in[0,1)$ be an irrational number and 
let $R_\alpha$ be the translation on $[0,1)$ defined 
by $R_\alpha(t)=t+\alpha$ ($\md 1$). 
We define the map $\pi:[0,1)\rightarrow \{0,1\}$ by 
$\pi(t)=0$ if $t\in[0,\alpha)$ and 
$\pi(t)=1$ otherwise. 
Let $X$ be the closure of 
$\{(\pi(R_\alpha^n(t)))_{n\in \Z}:t\in[0,1)\}$ in 
$\{0,1\}^\Z$ and 
let $\phi$ be the two-sided shift map on $X$. 
The dynamical system $\xp$ is called a Sturmian shift 
arising from the $\alpha$-rotation $R_\alpha$ on $[0,1)$, 
and known to be a minimal subshift. 
This class of shifts was introduced by Morse and Hedlund in 1940 
and has been studied by several authors. 
On this subject the reader can consult Section 13.7 of \cite{LM}. 

We will borrow the notation from the proof of Theorem \ref{Dfg}. 
Let $U=[0.]$ and 
let $\mu$ be a unique $\phi$-invariant probability 
measure on $X$. Then $\mu(U)=\alpha$. 
Since $\xp$ is a kind of Denjoy system, we know that 
its $K^0$-group $K^0\xp$ is isomorphic to $\Z\oplus\Z$ 
by \cite[Theorem 5.3]{PSS}. 
The integration by $\mu$ gives the isomorphism from $K^0\xp$ to 
$\Z+\Z\alpha\subset\R$, 
and $[1_X]$ and $[1_U]$ are generators of $K^0\xp$. 
By Theorem \ref{CMsignature} and \ref{CMsimple}, 
the commutator subgroup $D([[\phi]]_0)$ is simple and 
$[[\phi]]_0/D([[\phi]]_0)$ is isomorphic to 
$K^0\xp/2K^0\xp\cong\Z_2\oplus\Z_2$. 
Furthermore, 
$[[\phi]]$ is finitely generated by Corollary \ref{DfgCOR}. 

We will verify that $[[\phi]]$ is generated by three elements. 
(The number of generators is at least three, because 
$[[\phi]]/D([[\phi]])\cong\Z\oplus\Z_2\oplus\Z_2$ is of rank three.) 
A Sturmian shift arising from the $(-\alpha)$-rotation is 
flip conjugate to $\xp$, 
and so we may assume that $\alpha$ is less than one half 
without loss of generality. 
Then $U$ and $\phi(U)$ are disjoint. 
Let $n$ be the maximum natural number satisfying $(n+1)\alpha<1$. 
If $0$ appears in a two-sided sequence $x\in X$, 
then at least $n$ 1's follow it. 
We denote consecutive $n$ 1's by $1^n$ simply. 
When we put 
\[ V=[0\overbrace{11\dots1}^n1.]=[01^n1.], \]
it is easily seen that $V$ and $\phi(V)$ are disjoint, and 
$\mu(V)=1-(n+1)\alpha$. 
Define $\sigma_U,\sigma_V\in[[\phi]]_0$ by 
\[ \sigma_U(x) \ = \ \begin{cases}
\phi(x) & x\in U \\
\phi^{-1}(x) & x\in \phi(U) \\
x & \mbox{otherwise} \end{cases}\]
and 
\[\sigma_V(x) \ = \ \begin{cases}
\phi(x) & x\in V \\
\phi^{-1}(x) & x\in \phi(V) \\
x & \mbox{otherwise}. \end{cases}\]
Then the two elements 
$\sgn(\sigma_U)=[1_U]+2K^0\xp$ and $\sgn(\sigma_V)=[1_V]
+2K^0\xp$ generate $K^0\xp/2K^0\xp\cong\Z_2\oplus\Z_2$. 
We claim that $[[\phi]]$ is generated by $\sigma_U$, $\sigma_V$ and 
$\phi$. 
Let us consider the case $n\geq 2$ at first. 
Taking the commutator of $\sigma_U$ and $\phi\sigma_U\phi^{-1}$, 
we have $\gamma_{\phi(U)}$. 
The equations 
\begin{gather*}
\phi^{n-1}\gamma_{\phi(U)}\phi^{1-n}=\gamma_{\phi^n(U)}
=\gamma_{[01^n.]}\\
\sigma_U\gamma_{[01^n.]}^{-1}\sigma_U\gamma_{[01^n.]}=
\gamma_{[01^n0.]}
\end{gather*}
and 
\[ \sigma_V\gamma_{[01^n.]}^{-1}\sigma_V\gamma_{[01^n.]}=
\gamma_{[01^n1.0]} \]
are easily checked. 
Conjugating $\gamma_{[01^n0.]}$ and $\gamma_{[01^n1.0]}$ 
by some power of $\phi$, 
we get $\gamma_{[01^n01^n.]}$ and $\gamma_{[01^n101^n.]}$. 
By repeating this procedure, for every cylinder set $W$, 
$\gamma_W$ will be obtained. 
Thereby $[[\phi]]$ is generated. 
In the case $n=1$ we can do the same argument, 
starting from $\sigma_U\sigma_V\sigma_U\sigma_V=\gamma_{[0110.]}$. 
\end{exm}

\flushleft{
\textit{e-mail; matui@math.s.chiba-u.ac.jp \\
Graduate School of Science and Technology,\\
Chiba University,\\
1-33 Yayoi-cho, Inage-ku,\\
Chiba 263-8522,\\
Japan. }}

\end{document}